\documentclass[11pt]{amsart}
\usepackage{amssymb, latexsym, tikz}

\theoremstyle{plain}
\newtheorem{theorem}{Theorem}
\newtheorem{corollary}{Corollary}
\newtheorem*{thm-cheb}{Theorem (Chebyshev)}

\newtheorem{proposition}{Proposition}

\newtheorem*{2'}{Theorem 2'}
\newtheorem*{3'}{Theorem 3'}

\theoremstyle{remark}

\newtheorem*{Remark 1}{Remark 1}
\newtheorem*{Remark 2}{Remark 2}
\newtheorem*{Remark 3}{Remark 3}
\newtheorem*{Remark 4}{Remark 4}

\numberwithin{equation}{section}

\begin{document}

\title[ The Infinite Limit of Separable Permutations]
 {The Infinite Limit of Separable Permutations}

\author{Ross G. Pinsky}


\address{Department of Mathematics\\
Technion---Israel Institute of Technology\\
Haifa, 32000\\ Israel}
\email{ pinsky@math.technion.ac.il}

\urladdr{http://www.math.technion.ac.il/~pinsky/}

\subjclass[2000]{60C05, 60B10, 05A05} \keywords{separable permutations, random permutation, Schr\"oder numbers }
\date{}

\begin{abstract}
Let $P_n^{\text{sep}}$ denote the uniform probability measure on the set of separable permutations in $S_n$.
Let $\mathbb{N}^*=\mathbb{N}\cup\{\infty\}$ with an appropriate metric and denote by
 $S(\mathbb{N},\mathbb{N}^*)$ the compact metric space consisting  of functions
$\sigma=\{\sigma_i\}_{ i=1}^\infty$ from
$\mathbb{N}$ to $\mathbb{N}^*$  which
are injections when restricted to $\sigma^{-1}(\mathbb{N})$\rm; that is, if $\sigma_i=\sigma_j$, $i\neq j$, then
$\sigma_i=\infty$.
Extending permutations $\sigma\in S_n$
by defining $\sigma_j=j$, for $j>n$, we have $S_n\subset
S(\mathbb{N},\mathbb{N}^*)$.
We show that $\{P_n^{\text{sep}}\}_{n=1}^\infty$  converges weakly on $S(\mathbb{N},\mathbb{N}^*)$ to a limiting distribution of regenerative type, which we
calculate explicitly.

\end{abstract}

\maketitle
\section{Introduction and Statement of Results}

Let $S_n$ denote the permutations of $[n]:=\{1,\cdots, n\}$. Given $\sigma\in S_k$ and $\tau\in S_l$,
the \it direct sum\rm\ of $\sigma$ and $\tau$ is the permutation in $S_{k+l}$ given by
$$
(\sigma\oplus\tau)(i)=\begin{cases}\sigma(i),\ i=1,\cdots, k;\\ \tau(i-k)+k,\ i=k+1,\cdots k+l,\end{cases}
$$
and the \it skew sum\rm\ $\sigma\ominus\tau$
is the permutation in $S_{k+l}$ given by
$$
(\sigma\ominus\tau)(i)=\begin{cases}\sigma(i)+l ,\ i=1,\cdots, k;\\ \tau(i-k),\ i=k+1,\cdots k+l.\end{cases}
$$
A permutation is \it indecomposable\rm\ if it cannot  be represented as the direct sum of two nonempty permutations and is \it skew  indecomposable \rm  if it cannot  be written as the skew sum
of two nonempty permutations.
A permutation is \it separable \rm if it can be obtained from the singleton permutation   by iterating direct sums and skew sums.
Equivalently, a permutation is separable if it can be successively decomposed and skew decomposed until all of the indecomposable and skew  indecomposable pieces of the permutation are singletons.
(For example, using one-line notation, consider the separable permutation  $\sigma=4352167$. It can be decomposed into $43521\oplus12$. Then 43521 can be skew decomposed into 213$\ominus$21 and 12 can be
decomposed into $1\oplus1$.
Now 213 can be decomposed into 21$\oplus$1 and  21 can be skew decomposed into 1$\ominus$1. Finally, again 21 can be skew decomposed into $1\ominus$1.)
It is well-known \cite{BBL}   that a permutation is separable if and only it avoids the patterns 2413 and 3142.

Let $\text{SEP}(n)$ denote the set of separable  permutations in $S_n$, and let $P^{\text{sep}}_n$ denote the uniform measure on  SEP$(n)$.
In this paper we investigate the weak limiting behavior of the probability measures $\{P^{\text{sep}}_n\}_{n=1}^\infty$ as $n\to\infty$.
In the limit we  obtain a probability measure not on the set of permutations of $\mathbb{N}$,
but on a more general structure that we now describe.

Let $\mathbb{N}^*=\mathbb{N}\cup \{\infty\}$ with the metric $d_{N^*}(i,j)=\sum_{k=i}^{j-1}2^{-k}$, for $1\le i<j\le \infty$.
Denote by $S(\mathbb{N},\mathbb{N}^*)$ the set of functions
$\sigma=\{\sigma_i\}_{ i=1}^\infty$ from
$\mathbb{N}$ to $\mathbb{N}^*$ \it which
are injections when restricted to $\sigma^{-1}(\mathbb{N})$\rm; that is, if $\sigma_i=\sigma_j$, $i\neq j$, then
$\sigma_i=\infty$.
The space $S(\mathbb{N},\mathbb{N}^*)$ can be identified with
the countably infinite product $\mathbb{N}^*\times\mathbb{N}^*\cdots$.
Since $\mathbb{N}^*$ is a compact metric space, it follows that $S(\mathbb{N},\mathbb{N}^*)$
is also a compact metric space with the metric $D(\sigma,\tau):=\sum_{i=1}^\infty \frac{d_{\mathbb{N}^*}(\sigma_i,\tau_i)}{2^{i}}$.
Any permutation $\sigma\in S_n$ may be identified  uniquely with an element of $S(\mathbb{N},\mathbb{N}^*)$ by defining
$\sigma_j=j$, for $j>n$.
Consequently, if $\mu_n$ is a probability measure on $S_n$, for each $n\in\mathbb{N}$, then
$\{\mu_n\}_{n=1}^\infty$ may be considered as a  sequence of probability measures on
the compact metric space $S(\mathbb{N},\mathbb{N}^*)$.

Let $P_n$ denote the uniform measure on $S_n$. It is easy to see that the sequence
$\{P_n\}_{n=1}^\infty$, considered as probability measures on $S(\mathbb{N},\mathbb{N}^*)$,
converges weakly to the trivial $\delta$-measure that places all its mass on the function
$\sigma\in S(\mathbb{N},\mathbb{N}^*)$ satisfying $\sigma(j)=\infty,\ j\in\mathbb{N}$.
However, when one works with the uniform measure on certain classes of pattern avoiding permutations, one
obtains non-trivial limits.
We will obtain in very explicit form the weak limit of the probability measures
$\{P^{\text{sep}}_n\}_{n=1}^\infty$ considered as  probability measures on  $S(\mathbb{N},\mathbb{N}^*)$.
The limiting probability measure will be  in the form of a regenerative concatenation.
Since the topology on $S(\mathbb{N},\mathbb{N}^*)$ is the product topology,
 the convergence of  $\{P^{\text{sep}}_n\}_{n=1}^\infty$  to a limiting
probability measure
allows one to understand for any fixed $m$, the asymptotic behavior of the statistics of $\sigma_1\cdots\sigma_m$, where
$\sigma\in \text{SEP}(n)$ is uniformly random and $n\to\infty$.
We note that \cite{P} considered similar types of limits for permutations avoiding a particular pattern of length three.
The remark after Corollary \ref{cor} points out a fundamental difference between the limiting behavior obtained here for
separable permutations and that obtained in \cite{P} for permutations avoiding a particular pattern  of length three.

With the exception of  the class of permutations avoiding a particular pattern of length three, the class of separable permutation is the most studied class of pattern avoiding permutations. The fact that these permutations can be completely decomposed
and the fact that they can be enumerated by a closed form generating function (see below) make them tractable.
The study of general pattern avoiding permutations goes back to Knuth's observation \cite{K} that a permutation is so-called stack sortable if and only if it 231-avoiding.  Similarly, the study of separable permutations goes back to
\cite{AN} where it was shown that these are precisely the permutations which are sortable by so-called pop stacks.
Separable permutations also arise in a variety of other applications, for example in bootstrap percolation \cite{SS}
and in connection to   polynomial interchanges where one studies the possible ways that the relative order of the values of a family of polynomials can be modified when crossing a common zero \cite{G}.

We now set the stage  in order to describe our convergence result.
A finite set of consecutive integers in $\mathbb{N}$ will be called a \it block\rm.
For $a,b\in\mathbb{R}$ with $a\le b$, we will denote the block $\{a,\cdots, b\}$ by $[a,b]$.
As has  already been mentioned, the limiting distribution of $\{P_n^{\text{sep}}\}_{n=1}^\infty$ that  we will obtain on $S(\mathbb{N},\mathbb{N}^*)$ has
  a regenerative structure.
In order to describe this regenerative structure,
we  need to consider permutations of blocks $I$. Denote the set of permutations of a block $I$ by $S_I$.
Thus,
for example, if $I=[3,5]$, then there are six permutation in $S_I$; namely $(3\ 4 \ 5), (3\ 5\ 4),(4\ 3\ 5), (4\ 5\ 3), (5\ 3 \ 4),(5\ 4 \ 3)$.
We will denote a generic permutation  of $I$ by $\sigma_I$.
We define a  separable permutation in  $S_I$ in the obvious way, analogous to our original definition of a separable permutation, or equivalently, as a permutation
 that avoids the patterns 2413 and 3142---see the definition of $\text{pat}(\sigma_I)$ in the paragraph
 after the remark following Corollary \ref{cor}. Similarly, we define indecomposable and skew indecomposable separable permutations in $S_I$.
For example, the permutations  $(4\ 2\ 3\ 6\ 5)$ and $(5\ 4\ 6\ 2\ 3)$ of the block $[2,6]$ are both separable,
the former one being skew indecomposable  and
  the latter one begin indecomposable.
We also  define $\infty^{(j)}$ to be the $j$-fold image of $\infty$: $\infty^{(j)}=\underbrace{(\infty\infty\cdots\infty)}_{j\ \text{times}}$,  $j\in\mathbb{N}$; we call this a block of infinities.
We  can now concatenate these simple permutation  of blocks with each other and with blocks of infinities to build   more complicated objects.
For example, if $I_1=[3,5]$ and $I_2=[20,23]$, and if the permutations $\sigma_{I_i}$, $i=1,2$,
are given by
$\sigma_{I_1}=(5\ 3 \ 4)$ and  $\sigma_{I_2}=(22\ 20\ 21\ 23)$,
then
$$
\infty^{(2)}*\sigma_{I_1}*\infty^{(1)}*\sigma_{I_2}:=(\infty, \infty,5,3,4, \infty,22, 20, 21, 23).
$$

Let $s_n=|\text{SEP}(n)|,\ n\ge1$, denote the number of separable permutations in $S_n$.
Let
$$
s(x)=\sum_{n=1}^\infty s_nx^n
$$
denote the generating function of $\{s_n\}_{n=1}^\infty$.
For a separable permutation, define the length of the first indecomposable block and the length of the first skew indecomposable block respectively by
\begin{equation}\label{Bdef}
\begin{aligned}
&|B_1^{+,n}|(\sigma)=\min\{j:\sigma([j])=[j]\}, \ \sigma\in \text{SEP}(n);\\
&|B_1^{-,n}|(\sigma)=\min\{k:\sigma([k])=[n]-[n-k]\}, \sigma\in \text{SEP}(n).
\end{aligned}
\end{equation}
Let $B_1^{+,n}(\sigma)$ denote the corresponding permutation  of
the block \newline $\big[1,|B_1^{+,n}(\sigma)|\big]$ (the first $|B_1^{+,n}(\sigma)|$ entries in $\sigma$),
and let
$B_1^{-,n}(\sigma)$ denote the corresponding permutation  of
the block $\big[n-|B_1^{-,n}(\sigma)|+1, n\big]$ (the first $|B_1^{-,n}(\sigma)|$ entries in $\sigma$).
Note that, by construction, the permutation
$B_1^{+,n}(\sigma)$  is indecomposable
and the permutation  $B_1^{-,n}(\sigma)$ is skew decomposable.

By the definition of separable permutations,
for each $\sigma\in \text{SEP}(n)$, with $n\ge2$, exactly one out of
$|B_1^{+,n}|(\sigma)$ and $|B_1^{-,n}|(\sigma)$ is equal to $n$,
and by symmetry,
\begin{equation}\label{sym}
|\{\sigma\in \text{SEP}(n):|B_1^{+,n}|(\sigma)=n|=|\{\sigma\in \text{SEP}(n):|B_1^{-,n}|(\sigma)=n|=
\frac12 s_n,\ n\ge2.
\end{equation}
That is, half of the permutations in SEP$(n)$, $n\ge2$,  are indecomposable and half are skew indecomposable.
Partitioning SEP$(n)$ by $\{|B_1^{+,n}|=j\}_{j=1}^n$ (or alternatively, by $\{|B_1^{-,n}|=j\}_{j=1}^n$), and
using  the concatenating structure of separable permutations, it  follows that
\begin{equation}\label{rec}
s_n=s_1s_{n-1}+\frac12\sum_{j=2}^{n-1}s_js_{n-j}+\frac12 s_n,\ n\ge2.
\end{equation}
From this it is easy to show that
\begin{equation}\label{genfunc}
s(x)=\frac12(1-x-\sqrt{x^2-6x+1}\thinspace),\ \text{for}\ |x|<3-2\sqrt{2}.
\end{equation}
From the above formula for  the generating function,  one can show  that
\begin{equation}\label{asymp}
s_n\sim\frac1{2\sqrt{\pi n^3}}(3-2\sqrt2)^{-n+\frac12}.
\end{equation}
(See \cite[p. 474-475]{FS}. Our  $s_n$ is equal to their $D_{n-1}$, and $D_n$ is known at the $n$th Schr\"oder number (or ``big'' Schr\"oder number.)
From \eqref{asymp}, it follows that \eqref{genfunc} also holds for $|x|=3-2\sqrt2$; we have
\begin{equation}\label{svalue}
s(3-2\sqrt2)=\sqrt2-1.
\end{equation}

In \eqref{chi}-\eqref{unifsepindecom} below, we
 define  the five types of  random variables that will be used to describe the limiting distribution of  $\{P^{\text{sep}}_n\}_{n=1}^\infty$:
\begin{equation}\label{chi}
\begin{aligned}
&\{\chi_{0,1}^{(n)}\}_{n=1}^\infty \ \text{are IID and distributed according to the Bernoulli distribution}\\
&\text{with parameter}\ \frac12: \\
&P(\chi_{0,1}^{(n)}=0)=P(\chi_{0,1}^{(n)}=1)=\frac12;\\
\end{aligned}
\end{equation}
\begin{equation}\label{N}
\begin{aligned}
&\{N^{(n)}\}_{n=1}^\infty \ \text{are IID and distributed according to the following distribution:}\\
&  P(N^{(n)}=j)=\begin{cases}\frac{\sqrt2}2,\ j=0;\\
\sqrt2(3-2\sqrt2)^j,\ j=1,2,\cdots.\end{cases}
\end{aligned}
\end{equation}
\bf\noindent Remark.\rm\ For later use in Proposition \ref{stable},
we note that $EN^{(n)}=(\frac12)^\frac32$.

\begin{equation}\label{R}
\begin{aligned}
&\big\{\{R^{(n)}_m\}_{m=1}^\infty\big\}_{n=1}^\infty\ \text{are IID and distributed according to the following distribution:}\\
&P(R^{(n)}_m=k)=\frac{s_k(3-2\sqrt2)^k}{\sqrt2-1},\ k=1,2,\cdots.
\end{aligned}
\end{equation}
\begin{equation}\label{L}
\begin{aligned}
&\{L^{(n)}\}_{n=1}^\infty \ \text{are IID and distributed according to the following distribution:}\\
&P(L^{(n)}=j)=\begin{cases}\frac{2(3-2\sqrt2)}{2-\sqrt2};\ j=1\\ \frac{s_j(3-2\sqrt2)^j}{2-\sqrt2},\  j=2,3,\cdots;\end{cases}.
\end{aligned}
\end{equation}
\noindent \bf Remark.\rm\ It follows from \eqref{svalue} and the fact that $s_1=1$ that the distributions in \eqref{R} and \eqref{L} are indeed probability distributions.

\begin{equation}\label{unifsepindecom}
\begin{aligned}
&\Pi^{\text{sep}}(I)\ \text{is a uniformly random, indecomposable, separable permutation}\\
&\text{ of the finite block} \ I\subset \mathbb{N},
\text{and}\ \{\Pi^{\text{sep}}(I): I\subset\mathbb{N}\ \text{a finite block}\}\ \text{are independent}.
\end{aligned}
\end{equation}

All the random variables in \eqref{chi}-\eqref{unifsepindecom} are assumed to be mutually independent.

In the theorem below, we present a  rather involved  formula for the  $S(\mathbb{N},\mathbb{N}^*)$-valued random variable
 whose distribution is the limiting distribution of
$\{P^{\text{sep}}_n\}_{n=1}^\infty$.
It is worthwhile to begin with a more verbal and informal  description.

The random variable is formed by regenerative concatenation.
Its first piece is constructed via the random variables
$\chi_{0,1}^{(1)}, N^{(1)}, \{R_m^{(1)}\}_{m=1}^\infty, L^{(1)}$ and the random variable $\Pi^{\text{sep}}(I_1)$ with random block $I_1$
as specified below.
To construct this piece, first we use the random variables $N^{(1)}$ and $\{R^{(1)}_m\}_{m=1}^\infty$ and  discard the block $[1,\sum_{m=1}^{N^{(1)}}R^{(1)}_m]$, by which we mean that these numbers will not appear anywhere in the range of the
$S(\mathbb{N},\mathbb{N}^*)$-valued random variable. (Note though, that it is possible for this block to be empty since $N^{(1)}$ can be equal to 0.)
Then with the addition of the random variables $\chi_{0,1}^{(1)}$ and $L^{(1)}$, we build a block as follows. If $\chi^{(1)}_{0,1}=0$, then we lay down the block
$\infty^{(L^{(1)})}$, that is, a row of infinities of length $L^{(1)}$.
On the other hand, if $\chi_{0,1}=1$, then we lay down the block $\Pi^{\text{sep}}(I_1)$, a uniformly random, indecomposable, separable
permutation  of the random block $I_1:=[1+\sum_{m=1}^{N^{(1)}}R^{(1)}_m,L^{(1)}+\sum_{m=1}^{N^{(1)}}R^{(1)}_m]$ of length $L^{(1)}$.
This completes the construction of the first piece.

The second piece is constructed using the random variables
$\chi_{0,1}^{(2)}$, $N^{(2)}$,\newline $\{R_m^{(2)}\}_{m=1}^\infty$, $L^{(2)}$ and the random variable $\Pi^{\text{sep}}(I_2)$ with random block $I_2$ as specified below.
Note that from the construction of the  first piece, the block of numbers $[1,\sum_{m=1}^{N^{(1)}}R^{(1)}_m]$  was discarded.  Furthermore, if $\chi_{0,1}^{(1)}$ was equal to 1, then the block of numbers
$I_1=[1+\sum_{m=1}^{N^{(1)}}R^{(1)}_m,L^{(1)}+\sum_{m=1}^{N^{(1)}}R^{(1)}_m]$  was used in the construction of the first piece, in which case these numbers are also no longer available.
We start the construction of the second piece by
 discarding a block of length $\sum_{m=1}^{N^{(2)}}R^{(2)}_m$, starting with the smallest number available.
This discarded block is given by
$$
\big[\chi_{0,1}^{(1)}L^{(1)}+\sum_{m=1}^{N^{(1)}}R^{(1)}_m+1,\thinspace\chi_{0,1}^{(1)}L^{(1)}+\sum_{m=1}^{N^{(1)}}R^{(1)}_m+\sum_{m=1}^{N^{(2)}}R^{(2)}_m\big];
$$
these numbers will not appear anywhere in the range of the $S(\mathbb{N},\mathbb{N}^*)$-valued random variable.
Then with the addition of the random variables $\chi_{0,1}^{(2)}$ and $L^{(2)}$, we build a block as follows. If $\chi^{(2)}_{0,1}=0$, then we lay down the block
$\infty^{(L^{(2)})}$, that is, a row of infinities of length $L^{(2)}$.
However, if $\chi^{(2)}_{0,1}=1$, then we lay down the block $\Pi^{\text{sep}}(I_2)$,
a uniformly random, indecomposable, separable
permutation  of the random block
$$
I_2:=\big[1+\chi_{0,1}^{(1)}L^{(1)}+\sum_{m=1}^{N^{(1)}}R^{(1)}_m+\sum_{m=1}^{N^{(2)}}R^{(2)}_m,\thinspace L^{(2)}+\chi_{0,1}^{(1)}L^{(1)}+\sum_{m=1}^{N^{(1)}}R^{(1)}_m+\sum_{m=1}^{N^{(2)}}R^{(2)}_m\big]
$$
of length $L^{(2)}$. This completes the construction of the second piece.
The construction of the  random variable continues in this regenerative fashion.

We now state the theorem.
For convenience, we set
$$
\chi_{0,1}^{(0)}=L^{(0)}=0.
$$
\newpage
\begin{theorem}\label{thm1}
The distributions $\{P^{\text{sep}}_n\}_{n=1}^\infty$, considered on the space $S(\mathbb{N},\mathbb{N}^*)$, converge weakly to the distribution of the random variable

\begin{equation}\label{theformula}
\begin{aligned}
&*_{n=1}^\infty\big(\chi_{0,1}^{(n)}\thinspace\Pi^{\text{sep}}(I_n)+(1-\chi_{0,1}^{(n)})\thinspace\infty^{(L^{(n)})}\big):=\\
&\big(\chi_{0,1}^{(1)}\thinspace \Pi^{\text{sep}}(I_1)+(1-\chi_{0,1}^{(1)})\thinspace\infty^{(L^{(1)})}\big)*\big(\chi_{0,1}^{(2)}\thinspace\Pi^{\text{sep}}(I_2)+(1-\chi_{0,1}^{(2)})\thinspace\infty^{(L^{(2)})}\big)*\cdots,
\end{aligned}
\end{equation}
with
$$
I_n=\big[1+\sum_{k=1}^n\chi_{0,1}^{(k-1)}L^{(k-1)}+\sum_{k=1}^n\sum_{m=1}^{N^{(k)}}R^{(k)}_m,\thinspace   \chi_{0,1}^{(n)}L^{(n)}+\sum_{k=1}^n\chi_{0,1}^{(k-1)}L^{(k-1)}+\sum_{k=1}^n\sum_{m=1}^{N^{(k)}}R^{(k)}_m\big],
$$
where the random variables involved here are as in \eqref{chi}-\eqref{unifsepindecom}.
\end{theorem}
From the theorem and from the explanation of the construction of the random variable appearing in the theorem, note that in the first $n$ pieces of the
concatenation defining  the  random variable,  the numbers from 1 up to
$\sum_{k=1}^n\chi_{0,1}^{(k)}L^{(k)}+\sum_{k=1}^n\sum_{m=1}^{N^{(k)}}R^{(k)}_m$ have been
involved. Of these numbers, $\sum_{k=1}^n\sum_{m=1}^{N^{(k)}}R^{(k)}_m$ of them have been discarded and don't appear in the  random variable, and
$\sum_{k=1}^n\chi_{0,1}^{(k)}L^{(k)}$ of them do appear in the  random variable.
Thus, the ratio of the quantity of  numbers appearing to the quantity of  numbers appearing  or discarded in the first  $n$
pieces of the concatenation is
given by
\begin{equation}\label{ratio1}
\frac{\sum_{k=1}^n\chi_{0,1}^{(k)}L^{(k)}}{\sum_{k=1}^n\chi_{0,1}^{(k)}L^{(k)}+\sum_{k=1}^n\sum_{m=1}^{N^{(k)}}R^{(k)}_m}.
\end{equation}

Similarly, from the theorem it follows that in the first $n$ pieces of the concatenation defining the random variable,
the number of  integers appearing is $\sum_{k=1}^n\chi_{0,1}^{(k)}L^{(k)}$
and the number of infinities appearing is $\sum_{k=1}^n(1-\chi_{0,1}^{(k)})L^{(k)}$.
Thus the ratio of the number of integers appearing to the number of integers and infinities appearing in the
first $n$ pieces of the concatenation is given by
\begin{equation}\label{ratio2}
\frac{\sum_{k=1}^n\chi_{0,1}^{(k)}L^{(k)}}{\sum_{k=1}^n\chi_{0,1}^{(k)}L^{(k)}+\sum_{k=1}^n(1-\chi_{0,1}^{(k)})L^{(k)}}.
\end{equation}

We  consider the limiting behavior of the  ratios in \eqref{ratio1} and \eqref{ratio2}.
Note from \eqref{LR} and \eqref{asymp} that the random variables $L^{(k)}$ and $R^{(k)}_m$ do not have a finite first moment, so the law of large numbers does not hold for them.
Actually, they are in the domain of attraction of stable laws with parameter $\frac12$. We have the following proposition.
\bigskip

\begin{proposition}\label{stable}
\noindent i.
$$
\lim_{n\to\infty}\frac1{n^2}\sum_{k=1}^nL^{(k)}\stackrel{\text{dist}}{=}Z_L,
$$
where $Z_L$ is the one-sided stable distribution  with stability parameter $\frac12$ and with characteristic function
$$
\phi_{Z_L}(t)=E\exp(-itZ_L)=\exp\Big(-(\frac12)^\frac34|t|^\frac12\big(1+i\thinspace\text{\rm sgn}(t)\big)\Big).
$$
\noindent  ii.
$$
\lim_{n\to\infty}\frac1{n^2}\sum_{k=1}^nR^{(k)}_m\stackrel{\text{dist}}{=}Z_R,
$$
where $Z_R$ is the one-sided stable distribution  with stability parameter $\frac12$ and with characteristic function
$$
\phi_{Z_R}(t)=E\exp(-itZ_R)=\exp\Big(-(\frac12)^\frac14|t|^\frac12\big(1+i\thinspace\text{\rm sgn}(t)\big)\Big).
$$
\noindent iii.
$$
\lim_{n\to\infty}\frac1{n^2}\sum_{k=1}^n\chi_{0,1}^{(k))}L^{(k)}\stackrel{\text{dist}}{=}\lim_{n\to\infty}\frac1{n^2}\sum_{k=1}^n\sum_{m=1}^{N^{(k)}}R^{(k)}_m=Z,
$$
where $Z$ is the one-sided stable distribution  with stability parameter $\frac12$ and with characteristic function
$$
\phi_Z(t)=E\exp(-itZ)=\exp\Big(-(\frac12)^\frac74|t|^\frac12\big(1+i\thinspace\text{\rm sgn}(t)\big)\Big).
$$

\end{proposition}
Proposition \ref{stable} immediately gives  the   following corollary concerning the ratio in \eqref{ratio1}.
\begin{corollary}\label{cor}
The ratio of the quantity of  numbers appearing to the quantity of  numbers appearing or discarded in the first  $n$
pieces of the concatenation defining  the random variable in Theorem \ref{thm1} satisfies
$$
\lim_{n\to\infty}\frac{\sum_{k=1}^n\chi_{0,1}^{(k)}L^{(k)}}{\sum_{k=1}^n\chi_{0,1}^{(k)}L^{(k)}+\sum_{k=1}^n\sum_{m=1}^{N^{(k)}}R^{(k)}_m}\stackrel{\text{dist}}{=}\frac{Z_1}{Z_1+Z_2},
$$
where $Z_1$ and $Z_2$ are IID random variables distributed as $Z$ in part (iii) of Proposition \ref{stable}.
\end{corollary}

\bf\noindent Remark.\rm\ Corollary \ref{cor} presents a fundamental
difference between the structure of the limiting probability measure on
$S(\mathbb{N},\mathbb{N}^*)$  for separable permutations and the corresponding limiting behavior
studied in \cite{P} in the case of permutations avoiding a pattern of length three.
The corollary above indicates
that the ratio of the quantity of numbers appearing  to the  quantity of numbers either appearing or discarded
in the first $n$ pieces of the concatenation converges to a limiting \it nondeterministic\rm\ quantity.
In \cite{P},
 for  three out of the six permutations in $S_3$, namely   312, 231 and 213, the limiting probability measure has
 a regenerative concatenation structure as it does here. It is easy to check
 from the structures obtained  there that
 the ratio of the quantity of numbers appearing  to the  quantity of numbers either appearing or discarded
in the first $n$ pieces of the concatenation converges in distribution to 1 as $n\to\infty$ in the case of 312,  is identically equal to
1 for all $n$ in the case of 231, and converges in distribution  to 0 as $n\to\infty$ in the case of 213.

\medskip
The next corollary shows that the limiting distribution of the ratio in \eqref{ratio2} is the same as that for \eqref{ratio1}.
We will use the proof of Proposition \ref{stable} to prove this.
\begin{corollary}\label{cortoprove}
The ratio of the number of integers appearing to the number of integers or infinities appearing in the first $n$ pieces
of the concatenation defining the random variable in Theorem \ref{thm1} satisfies
$$
\lim_{n\to\infty}\frac{\sum_{k=1}^n\chi_{0,1}^{(k)}L^{(k)}}
{\sum_{k=1}^n\chi_{0,1}^{(k)}L^{(k)}+\sum_{k=1}^n(1-\chi_{0,1}^{(k)})L^{(k)}}=\frac{ Z_1}{ Z_1+ Z_2},
$$
where $Z_1$ and $Z_2$ are IID   random variables distributed
as $Z$ in part (iii) of Proposition \ref{stable}.
\end{corollary}
\noindent \bf Remark.\rm\
Among other things, the two corollaries above express certain aspects of the symmetries inherent in separable permutations.
Recall that the reverse of a permutation $\sigma=\sigma_1\cdots\sigma_n\in S_n$, denoted by  $\sigma^{\text{rev}}$,  satisfies $\sigma^{\text{rev}}_j=\sigma_{n+1-j}$, $j\in[n]$.
Recall that the complement of $\sigma$, denoted
by $\sigma^{\text{comp}}$,  satisfies $\sigma^{\text{comp}}_j=n+1-\sigma_j$, $j\in[n]$. 
Finally, the reverse-complement of $\sigma$, denoted by $\sigma^{\text{rev-comp}}$,
 satisfies $\sigma^{\text{rev-comp}}_j=n+1-\sigma_{n+1-j}$, $j\in[n]$.
It is clear from the definition of a separable permutation that if $\sigma$ is distributed 
according to  the uniform measure $P_n^{\text{sep}}$ on the set SEP$(n)$ of separable permutations in $S_n$,
then  the three random variables $\sigma^{\text{rev}},\sigma^{\text{comp}}$ and $\sigma^{\text{rev-comp}}$ all have the same
distribution as the random variable $\sigma$. Furthermore, a weak convergence result similar to that in Theorem \ref{thm1} for
$\sigma$ can be given
for the four-tuple $(\sigma,\sigma^{\text{rev}},\sigma^{\text{comp}},\sigma^{\text{rev-comp}})$.
Of course all of the one-dimensional marginal distributions of the limiting random vector will coincide with the limiting distribution in
Theorem \ref{thm1}. 
This point of view can be used to gain some insight with regard to the numbers appearing, the numbers discarded and the infinities appearing in the limiting object in Theorem \ref{thm1} and with regard to the symmetry property
of the limiting random variable appearing in the two corollaries above. In order to describe this, we need to refer to the proof of Theorem \ref{thm1}. Thus, this discussion is postponed and  appears in the remark at the end 
of section \ref{thmpf}, after the proof of Theorem \ref{thm1}.

\medskip

 We end this first section by noting a type of limiting result, completely different from the type considered in Theorem \ref{thm1},
 that has been studied for the class of separable permutations as well as for the class of
permutations avoiding a particular pattern of length three.
 Let $\tilde S_n$ denote a class of permutations in $S_n$ as above, that is, the class of permutations that avoid a particular pattern  in $S_3$, or alternatively, that are separable.
Fix $m\ge 2$. For $\sigma\in \tilde S_n$, with $n>m$, and $I=\{i_1,i_2,\cdots, i_m\}$ with $1\le i_1<\cdots<i_m\le n$, let $\sigma_I=\sigma_{i_1}\cdots\sigma_{i_m}$.
Let $\text{pat}(\sigma_I)$ denote the permutation in $S_m$ which describes  the pattern of $\sigma_I$.
(For example, if $m=4, n=5$, $\sigma=32541$ and $I=\{1,3,4,5\}$, then
$\text{pat}(\sigma_I)=2431\in S_4$.)
Now fix a permutation $\pi\in S_m$. Let
$
\text{occ}_\pi(\sigma)=|\{I\subset [n] \ \text{of cardinality }\ m\ \text{such that}\ \text{pat}(\sigma_I)=\pi\}\}|$
denote  the number of occurrences of the pattern $\pi$ from among
the $\binom{n}{m}$ different $\sigma_I$'s of length $|I|=m$.
Considering $\sigma$ to be a uniformly random element of $\tilde S_n$, a number of papers have studied the asymptotic
behavior of $\text{occ}_\pi(\sigma)$ as $n\to\infty$.
In the case that $\tilde S_n$ is the class of permutations in $S_n$ that avoid a particular pattern of length three, this has been done at
  the level of expected value in a number of papers, see for example, \cite{B10},\cite{B12},\cite{H}, \cite{R}, while
  at the level of distribution this has been carried out in \cite{J}.
In the case that $\tilde S_n$ is the class of separable permutations, this has been done at the level of distribution
in \cite{BBFGP}.

The type of asymptotic behavior studied in this paper and the type of asymptotic behavior described in the  paragraph above
are mutually exclusive and complement one another. Indeed, as already noted, the type of asymptotic behavior studied in
this paper allows one to understand for any fixed $m$, the statistics of $\sigma_1\cdots\sigma_m$, where
$\sigma\in \text{SEP}(n)$ is uniformly random and $n\to\infty$.
However, the  asymptotic results
 described in the above paragraph give no information about this statistic.
Conversely, the behavior of   $\sigma_1\cdots\sigma_m$, for fixed $m$
has no influence on  the asymptotic results
 described in the above paragraph.

In section \ref{prelim} we proof some preliminary results that will be used in the proof of Theorem \ref{thm1}. In section \ref{thmpf} we prove Theorem \ref{thm1}
and then remark on how the symmetries concerning $\sigma, \sigma^{\text{rev}},\sigma^{\text{comp}}$
and $\sigma^{\text{rev-comp}}$ shed light on certain aspects of the theorem.
  In section \ref{proppf} we prove Proposition \ref{stable}
and Corollary \ref{cortoprove}.
\section{Preliminary results concerning separable permutations and Schr\"oder numbers}\label{prelim}

From \eqref{sym}, we have
\begin{equation}\label{symmetry}
P_n^{\text{SEP}}(|B_1^{+,n}|=n)=P_n^{\text{SEP}}(|B_1^{-,n}|=n)=\frac12,\ n\ge2,
\end{equation}
and from \eqref{rec} and the text immediately preceding it, we have
\begin{equation}\label{limiting+}
P_n(|B^{+,n}_1|=j\ \big|\thinspace|B^{+,n}_1|<n)=P_n(|B^{-,n}_1|=j\ \big|\thinspace|B^{-,n}_1|<n)=\begin{cases}2\frac{s_{n-1}}{s_n}, j=1;\\ \frac{s_js_{n-j}}{s_n},\ j=2,\cdots, n-1.\end{cases}
\end{equation}
From \eqref{symmetry} and \eqref{limiting+}, it follows
that the distributions of $|B_1^{+,n}|$ and $|B_1^{-,n}|$ under $P_n^{\text{SEP}}$ are identical; thus, in the sequel we will sometimes use the notation
$|B_1^{\pm,n}|$ in formulas that hold for both $|B_1^{+,n}|$ and $|B_1^{-,n}|$.
It also follows from \eqref{symmetry} and \eqref{rec} that
\begin{equation}\label{limiting-}
P_n(n-|B^{-,n}_1|=k\big|\thinspace |B^{-,n}_1|<n)= \begin{cases}\frac{s_ks_{n-k}}{s_n},\ k=1,2,\cdots n-2;\\ 2\frac{s_{n-1}}{s_n},\  k=n-1.\end{cases}
\end{equation}
(Formula \eqref{limiting-}  also holds with $|B^{+,n}_1|$ in place of $|B^{-,n}_1|$, but it won't be needed.)


From \eqref{asymp}, \eqref{limiting+} and \eqref{limiting-}, we have
\begin{equation}\label{limit+}
\lim_{n\to\infty}P_n(|B^{\pm,n}_1|=j\big|\thinspace |B^{\pm,n}_1|<n)=
\begin{cases}2(3-2\sqrt2),\ j=1;\\ s_j(3-2\sqrt2)^j,\  j=2,3,\cdots\end{cases}.
\end{equation}
and
\begin{equation}\label{limit-}
\lim_{n\to\infty}P_n(n-|B^{-,n}_1|=j\big|\thinspace |B^{-,n}_1|<n)=s_k(3-2\sqrt2)^k,\ k=1,2,\cdots.
\end{equation}
Using \eqref{svalue},
the sum on the right hand side of \eqref{limit+} is given by
\begin{equation}\label{sum+}
2(3-2\sqrt2)+\sum_{j=2}^\infty s_j(3-2\sqrt2)^j=3-2\sqrt2+s(3-2\sqrt2)=2-\sqrt2,
\end{equation}
and the sum on the right hand side of \eqref{limit-} is given by
\begin{equation}\label{sum-}
\sum_{k=1}^\infty s_k(3-2\sqrt2)^k=s(3-2\sqrt2)=\sqrt2-1.
\end{equation}
From \eqref{limit+} and \eqref{sum+} it follows that the random variables  $\{|B^{\pm,n}_1|\}_{n=1}^\infty$ under the measures
$P_n(\cdot\ \big| \thinspace |B^{\pm,n}_1|<n)$, considered on the space $\mathbb{N}^*$, converge weakly to
the random variable $\tilde L$ with distribution
\begin{equation}\label{tildeL}
P\big(\tilde L=j)=\begin{cases}2(3-2\sqrt2);\ j=1;\\ s_j(3-2\sqrt2)^j,\  j=2,3,\cdots;\\ \sqrt2-1,\ j=\infty.\end{cases}
\end{equation}
From \eqref{limit-} and \eqref{sum-} it follows that the random variables  $\{n-|B^{-,n}_1|\}_{n=1}^\infty$
under the measures
$P_n(\cdot\ \big| \thinspace |B^{-,n}_1|<n)$, considered on the space $\mathbb{N}^*$, converge weakly to
the random variable $\tilde R$ with distribution
\begin{equation}\label{tildeR}
P\big(\tilde R=j)=\begin{cases}\frac{s_k(3-2\sqrt2)^k}{\sqrt2-1}\, k=1,2,\cdots;\\ 2-\sqrt2,\ j=\infty.\end{cases}
\end{equation}

Since  $(2-\sqrt2)+(\sqrt2-1)=1$, it follows that
 on the space $\mathbb{N}^*\times\mathbb{N}^*$, the  random vectors $\{(|B_1^{-,n}|,n-|B_1^{-,n}|)\}_{n=1}^\infty$
under the measures $P_n(\cdot\ \big| \thinspace |B^{-,n}_1|<n)$ converge to the random vector $(\tilde L,\tilde R)$ with distribution
\begin{equation}\label{2dimdist}
\begin{aligned}
&P\big((\tilde L,\tilde R)=(j,\infty)\big)=\begin{cases}2(3-2\sqrt2),\ j=1;\\ s_j(3-2\sqrt2)^j,\  j=2,3,\cdots;\end{cases}\\
&P\big( (\tilde L,\tilde R)=(\infty,k)\big)=s_k(3-2\sqrt2)^k,\ k=1,2,\cdots.
\end{aligned}
\end{equation}
Note in particular that $P\big( (\tilde L,\tilde R)=(\infty,\infty)\big)=0$.

Denote respectively  by $R$ and $L$ random variables whose distributions are those of $R$ conditioned on $\{\tilde R<\infty\}$
and $L$ conditioned on $\{\tilde L<\infty\}$. That is,
\begin{equation}\label{LR}
\begin{aligned}
&P( R=k)=\frac{s_k(3-2\sqrt2)^k}{\sqrt2-1}\, k=1,2,\cdots;\\
&P( L=j)=\begin{cases}\frac{2(3-2\sqrt2)}{2-\sqrt2};\ j=1;\\ \frac{s_j(3-2\sqrt2)^j}{2-\sqrt2},\  j=2,3,\cdots.\end{cases}
\end{aligned}
\end{equation}
Note that $R$ and $L$ are the distributions respectively  of $R^{(n)}_m$ and  $L^{(n)}$  in \eqref{R} and \eqref{L}.

\section{Proof of Theorem \ref{thm1} }\label{thmpf}

To write down a complete and entirely rigorous proof of the theorem is extremely tedious and may well obscure the relative simplicity
of the ideas behind the proof. Thus, we will give a somewhat informal proof, with quite a bit of verbal explanation, relating at times simultaneously
to the situation for large $n$ and the situation in the limit after  $n\to\infty$.
At the end of this proof, we
then   prove completely rigorously a particular case  of the proof. From this, it will be clear that
one can precede similarly to obtain the entire proof.

Consider a permutation $\sigma$ under $P_n^{\text{sep}}$, for $n$ very large. That is, $\sigma$ is a uniformly distributed
separable permutation in $S_n$.
By \eqref{symmetry}, with probability $\frac12$, one has $|B_1^{+,n}|(\sigma)<n$ and $|B_1^{-,n}|(\sigma)=n$, and with probability $\frac12$ one has
$|B_1^{+,n}|(\sigma)=n$ and $|B_1^{-,n}|(\sigma)<n$.
Consider first the former case. Then the distribution of $|B_1^{+,n}|(\sigma)$ (under $P_n(\cdot\ \big| \thinspace |B^{+,n}_1|<n)$)   is given by
the distribution in \eqref{limiting+}.
From \eqref{limit+} and \eqref{sum+}, it follows that as $n\to\infty$,
$|B_1^{+,n}|(\sigma)$ will run off to infinity with probability approaching $\sqrt2-1$, while with probability approaching $2-\sqrt2$, it will converge to a
limit which is distributed as $L$ in  \eqref{LR}.
Now consider the latter case. The distribution of $|B_1^{-,n}|(\sigma)$ (under $P_n(\cdot\ \big| \thinspace |B^{-,n}_1|<n)$)  is also given by
the distribution in \eqref{limiting+}.
Thus, as with $|B_1^{+,n}|(\sigma)$,
$|B_1^{-,n}|(\sigma)$ will run off to infinity with probability approaching $\sqrt2-1$, while with probability approaching $2-\sqrt2$, it will converge to a
limit which is distributed as  $L$ in  \eqref{LR}.

We denote the four mutually exclusive cases above as follows:

\begin{equation}\label{4cases}
\begin{aligned}
& (\bf+F\rm):\ \  |B_1^{+,n}(\sigma)|<n\ \text{and}\ |B_1^{+,n}(\sigma)|\ \text{does not run off to infinity;}\\
& (\bf-F\rm):\ \  |B_1^{-,n}(\sigma)|<n \ \text{and}\ |B_1^{-,n}(\sigma)|\ \text{does not run off to infinity;}\\
& (\bf+I\rm):\ \  |B_1^{+,n}(\sigma)|<n\ \text{and} \ |B_1^{+,n}(\sigma)|\ \text{runs off to infinity;}\\
&(\bf -I\rm):\ \ |B_1^{-,n}(\sigma)|<n\ \text{and}\ |B_1^{-,n}(\sigma)|\ \text{runs off to infinity.}
\end{aligned}
\end{equation}
The probabilities for these four cases are respectively
\begin{equation}\label{4probs}
p_{(+F)}=\frac12(2-\sqrt2);\ \
  p_{(-F)}=\frac12(2-\sqrt2);\ \    p_{(+I)}=\frac12(\sqrt2-1);\ \   p_{(-I)}=\frac12(\sqrt2-1).
\end{equation}
The notation (+) will denote the union of the two cases (+F)
and (+I), and the notation $(-)$ will denote the union of the two cases ($-$F) and ($-$I).
Note that the cases (+) are the cases in which the  permutation $\sigma$ is skew indecomposable and the cases ($-$) are the cases in
which it is indecomposable. The notation (F) will denote the union of the two cases (+F) and ($-$F).

If (F) occurs,  then the first piece of the concatenation in the statement of the theorem
is obtained immediately.
Indeed, recalling the definition of
$B_1^{\pm,n}(\sigma)$ after \eqref{Bdef},
and recalling \eqref{limit+}, \eqref{LR} and \eqref{L},
we see that in the  case (+F), this piece is a uniformly random indecomposable separable permutation
$\Pi^\text{sep}(I_1)$, where $I_1=[1,L^{(1)}]$, while in the  case ($-$F) it is the block
of infinities $\infty^{(L^{(1)})}$, where $L^{(1)}$ is as in the statement of the theorem.
Setting this piece aside now, in the case (+F), this essentially returns us to the situation we started from, except that now we are looking
at uniformly random separable permutations of $[L^{(1)}+1,n]$, for large $n$, while in the  case ($-$F), this returns us exactly to
the situation we started from, and we again look at uniformly random separable permutations of $[n]$, for large $n$.
We emphasize that if (F) occurs, then we obtain the first piece of the concatenation and with regard to the rest
of the permutation, we are again presented with the four possible cases in \eqref{4cases} with the four corresponding probabilities in \eqref{4probs}.

Now consider the case (+I).
Then for large $n$,
$m:=|B_1^{+,n}(\sigma)|$  will be very large itself. Recalling the definition of
$B_1^{+,n}(\sigma)$ after \eqref{Bdef}, we see that $B_1^{+,n}(\sigma)$ is a uniformly distributed indecomposable separable permutation
in $S_m$ with $m$ large.
Considering \eqref{4cases} with $B_1^{+,n}(\sigma)$ in place of $\sigma$ and with large $m=|B_1^{+,n}(\sigma)|$ in place of large $n$,
this moves   us over to one of the two cases in ($-$), with corresponding probabilities
$2p_{(-F)}$ for moving to $(-F)$  and $2p_{(-I)}$ for moving to $(-I)$.
 Note that here the first piece of the concatenation
has not yet been obtained.

Now consider the case ($-$I).
Since we are assuming that  $|B_1^{-,n}(\sigma)|$ is running off to infinity,
 it follows  from the paragraph in which  \eqref{2dimdist}
appears that $n-|B_1^{-,n}(\sigma)|$  will converge to a    limit which is distributed as $R$ in \eqref{LR}.
By the definition   of $|B_1^{-,n}(\sigma)|$, the  block of numbers $\big[1,n-|B_1^{-,n}(\sigma)|\big]$ will appear in the final $n-|B_1^{-,n}(\sigma)|$ positions of $\sigma$. Thus, in the limit as $n\to\infty$, the first $R^{(1)}$ numbers in the permutation get swept out to infinity and will not appear in the limiting object, where $R^{(1)}$ is as in the
statement of the theorem.
Also, recalling the definition of
$B_1^{-,n}(\sigma)$ after \eqref{Bdef}, we see that $B_1^{-,n}(\sigma)$
is a uniformly distributed skew indecomposable  permutation  of
the block $\big[n-|B_1^{-,n}(\sigma)|+1, n\big]$. Since $n-|B_1^{-,n}(\sigma)$  converges weakly to $R^{(1)}$,
this essentially moves us over  to  one of the  two cases in (+), with corresponding probabilities
$2p_{(+F)}$ for moving to $(+F)$ and $2p_{(+I)}$ for moving to $(+I)$, except that now  we are looking at
uniformly random indecomposable separable    permutations of $[R^{(1)}+1,n]$, for large $n$.
Note that here the first piece of the concatenation
has not yet been obtained.

We summarize the mechanism we have discovered above in the following four statements:

\noindent \bf 1.\rm\ One begins from one of the four states in \eqref{4cases} with corresponding probabilities in \eqref{4probs}.
The first piece of the concatenation is obtained when a state in (F) is first reached.
If (+F) is reached, then the first part of the concatenation is an indecomposable separable permutation
whose length is distributed as $L$, and if ($-$F) is reached, then the first part of the concatenation is a block
of infinities
whose length is distributed as $L$.

\noindent \bf 2.\rm\ From the state (+I), one moves to the state ($-$F) with probability $2p_{-F}$, in which case the first piece of the concatenation
is obtained and it is a block
of infinities
whose length is distributed as $L$, while one moves to the state ($-$I) with probability $2p_{-I}$.

\noindent \bf 3.\rm\ From the state ($-$I), one moves to the state (+F) with probability $2p_{+F}$, in which case the first piece of the concatenation
is obtained and it is an indecomposable permutation whose length is distributed as $L$, while one moves to the state
(+I) with probability $2p_{+I}$.

\noindent \bf 4.\rm\  Every arrival at the  state ($-$I) causes the $r$ smallest numbers currently available (that is, that have not yet been involved in the construction) to be discarded; they will never be used in the construction. Here  $r$ has the distribution of $R$.
\medskip

Since $p_{+F}=p_{-F}$ and $p_{+I}=p_{-I}$, it is clear that with probability $\frac12$  the first piece of the concatenation will
be an indecomposable separable permutation, and  with probability $\frac12$ it will be a block of infinities. These two possibilities are represented in the statement of the theorem through the random variable
$\chi_{0,1}^{(1)}$.

We now show that the number of arrivals at state ($-$I) prior to the first  arrival at a state in (F) has the distribution
of $N^{(1)}$ as in \eqref{N}.
Let $A$ denote the number of such arrivals.
For $k\ge1$, in order to have $A\ge k$, either the first state visited is ($-$I) and then $k-1$ times in a row one
moves from ($-$I) to (+I) and back to ($-$I), or the first state visited is (+I), the next one is ($-$I) and then
$k-1$ times in a row one
moves from ($-$I) to (+I) and back to ($-$I).
The probability of the former scenario is $\frac12(\sqrt2-1)\big((\sqrt2-1)^2\big)^{k-1}$ while the probability
of the latter scenario is $\big(\frac12(\sqrt2-1)\big)^2\big((\sqrt2-1)^2\big)^{k-1}$.
Adding these, one obtains $P(A\ge k)=\frac{\sqrt2}2(\sqrt2-1)(\sqrt2-1)^{2k-2}$.
After a bit of arithmetic, one finds that
$$
P(A=k)=P(A\ge k)-P(A\ge k+1)=\sqrt2(3-2\sqrt2)^k, \ \text{for}\ k\ge1.
$$
Since $\sum_{k=1}^\infty\sqrt2(3-2\sqrt2)^k=1-\frac{\sqrt2}2$, one obtains $P(A=0)=\frac{\sqrt2}2$.
Thus the number of arrivals at state ($-$I) prior to the first arrival at a state in (F) indeed has the  distribution of $N^{(1)}$.

In light of the above, we see that by the time the first piece of the concatenation is constructed,
a certain block of  numbers, beginning from 1,  will have been removed from consideration, and the length of that block is distributed as $\sum_{m=1}^{N^{(1)}}R^{(1)}_m$.
Thus, in the case that the first piece of the concatenation is a permutation, it will be an indecomposable
separable permutation  of the random block $I_1=[1+\sum_{m=1}^{N^{(1)}}R^{(1)}_m,L^{(1)}+ \sum_{m=1}^{N^{(1)}}R^{(1)}_m]$.

After the first piece of the concatenation is obtained, be it a random permutation  or a block of infinities, we begin again with the  same mechanism,
however now the first number available for use is $1+\chi_{0,1}^{(1)}L^{(1)}+\sum_{m=1}^{N^{(1)}}R^{(1)}_m$.
In light of the  regenerative nature described above, this completes the proof of the theorem.

As promised at the beginning of the  proof,
we now give a completely rigorous derivation for one case, the case corresponding to arriving first at ($-$I) and then going to
(+F). Because the first state is assumed to be in ($-$),
 we are assuming that $|B_1^{-,n}|(\sigma)<n$. Let $n_1(\sigma)=|B_1^{-,n}|(\sigma)$.
Let $\tau=\tau(\sigma)$ denote the permutation  of the block
 $[n-n_1(\sigma)+1,n]$ corresponding to
the first $n_1(\sigma)$ entries of $\sigma$.
Note that $\tau$ is a skew indecomposable separable permutation. In fact,
 conditioned on $n-n_1(\sigma)=n-n_1$, $\tau$ is a uniformly distributed
skew indecomposable separable permutation of  the block  $[n-n_1+1,n]$.

We now take the liberty to extend in the obvious way
the domain of definition of $|B_1^{+,k}|(\cdot)$, $k\ge1$, which has been defined on separable  permutations in $S_k$, to include separable permutations of blocks of length $k$.
Since $\tau$ is a skew indecomposable  separable permutation, we have
$|B_1^{+,n_1(\sigma)}|(\tau)<n_1(\sigma)$.
Let $n_1'(\tau)=|B_1^{+,n_1(\sigma)}|(\tau)$, and let $\tau'=\tau'(\tau)$
denote the first $n_1'(\tau)$ entries of $\tau$. Note that $\tau'$ is an indecomposable separable permutation  of
the block $[n-n_1(\sigma)+1,n-n_1(\sigma)+n_1'(\tau)]$.
 In fact,
 conditioned on $n_1(\sigma)=n_1$ and $n_1'(\tau)=n_1'$, $\tau'$ is a uniformly distributed
indecomposable separable permutation of the block $[n-n_1+1, n-n_1+n_1']$.

For $j,k\ge1$ and for $n$ sufficiently large to accommodate the $j$ and $k$, we have
\begin{equation}\label{setup}
\begin{aligned}
&P_n^{\text{sep}}(n-n_1(\sigma)=k,\ n_1'(\tau)=j, \ \tau'=\tau_0)=\\
&P_n^{\text{sep}}(\tau'=\tau_0\big|\  n_1'(\tau)=j,\
n-n_1(\sigma)=k)\times\\
&P_n^{\text{sep}}(n_1'(\tau)=j|\  n-n_1(\sigma)=k)
\thinspace P_n^{\text{sep}}(n- n_1(\sigma)=k),\\
&\text{for}\ \tau_0\ \text{an indecomposable separable permutation  of}\ [k+1,k+j].
\end{aligned}
\end{equation}
Now
\begin{equation}\label{unifterm}
P_n^{\text{sep}}(\tau'=\tau_0\big|\  n_1'(\tau)=j,\
n-n_1(\sigma)=k)=\begin{cases} 1,\ j=1;\\ \frac2{s_j},\ j\ge2,\end{cases}
\end{equation}
 since there are $\frac{s_j}2$ indecomposable separable permutations of $[k+1,k+j]$, for $j\ge2$.
Also,
\begin{equation}\label{nextterm}
P_n^{\text{sep}}(n_1'(\tau)=j|\ n-n_1(\sigma)=k)=P_n^{\text{sep}}(|B_1^{+,n-k}|(\tau)=j\big|\ |B_1^{+,n-k}|(\tau)<n-k).
\end{equation}
(The conditioning on $\{|B_1^{+,n-k}|(\tau)<n-k\}$ comes in because $\tau$ is skew indecomposable.)
From \eqref{limit+} we have
\begin{equation}\label{nexttermagain}
\lim_{n-k\to\infty}P_n^{\text{sep}}(|B_1^{+,n-k}|(\tau)=j\big|\ |B_1^{+,n-k}|(\tau)<n-k)=\begin{cases}2(3-2\sqrt2),\ j=1;\\ s_j(3-2\sqrt2)^j,\  j=2,3,\cdots.\end{cases}
\end{equation}
Also, from \eqref{symmetry} and \eqref{limit-}, we have
\begin{equation}\label{last}
\lim_{n\to\infty}P_n^{\text{sep}}(n-n_1(\sigma)=k)=\lim_{n\to\infty}P_n^{\text{sep}}(n-|B_1^{-,n}|(\sigma)=k)=\frac12s_k(3-2\sqrt2)^k.
\end{equation}
From \eqref{setup}-\eqref{last}
we  obtain
\begin{equation}\label{threeprobs}
\begin{aligned}
&\lim_{n\to\infty}P_n^{\text{sep}}(n-n_1(\sigma)=k,\ n_1'(\tau)=j, \ \tau'=\tau_0)=\\
&\begin{cases}
\big(\frac12s_k(3-2\sqrt2)^k\big)\big(2(3-2\sqrt2)\big),\ j=1; k=1,2,\cdots;\\
\big(\frac12s_k(3-2\sqrt2)^k\big)\big(s_j(3-2\sqrt2)^j\big)\frac2{s_j}, \ j=2,3,\cdots; k=1,2,\cdots.\end{cases}
\end{aligned}
\end{equation}
We rewrite the right hand side of \eqref{threeprobs} as
\begin{equation}\label{finalrigor}
\begin{cases}
\big(\frac12(\sqrt2-1)\big)\big(\frac{s_k(3-2\sqrt2)^k}{\sqrt2-1}\big)  \big(2-\sqrt2\big)\frac{2(3-2\sqrt2)}{2-\sqrt2},\ j=1; k=1,2,\cdots;\\
\big(\frac12(\sqrt2-1)\big)(\frac{s_k(3-2\sqrt2)^k}{\sqrt2-1}\big)\big(2-\sqrt2\big)\frac{s_j(3-2\sqrt2)^j}{2-\sqrt2}\frac2{s_j}, \ j=2,3,\cdots; k=1,2,\cdots.\end{cases}
\end{equation}
The term $\frac12(\sqrt2-1)$ is
  $p_{(-I)}$  (corresponding to having started from $(-I)$), the term $2-\sqrt2$ is $2p_{(+F)}$ (corresponding to having moved from $(-I)$ to $(+F)$),
   the terms $\frac{s_k(3-2\sqrt2)^k}{\sqrt2-1}$, $k=1,2,\cdots$, give the distribution of $R$ as in \eqref{LR}, the terms
  $\frac{2(3-2\sqrt2)}{2-\sqrt2},\ j=1; \frac{s_j(3-2\sqrt2)^j}{2-\sqrt2},\ j=2,3,\cdots$
give the distribution of $L$ as in \eqref{LR}, and the term $\frac2{s_j},\ j\ge2$, indicates choosing uniformly
from the indecomposable separable permutations of the block  $[k+1,k+j]$.
    \hfill$\square$
    \medskip
    
\noindent \bf Remark.\rm\ We now return to discuss Theorem \ref{thm1}
in light of the symmetries noted in the remark following Corollary \ref{cortoprove}.
In the proof of the theorem, we delineated four cases in \eqref{4cases}, which were  denoted by $(+ F)$, $(-F)$, $(+I)$ and $(-I)$.
We represent them schematically in Figure 1.

From the explanation at the beginning of the proof of the theorem,
we see that  in the case of $(+F)$, the small box at the  lower left
indicates that a finite block of numbers will appear in the limiting object.
In the case of $(-F)$, the small box at the upper left indicates that a finite block of infinities will appear in the limiting object. In the case of $(+I)$,  the small box on the upper right indicates that this case has no effect on the limiting object. In the case of $(-I)$,  the small box on the lower right indicates that a finite block of numbers is discarded and will not appear in  the limiting object.
After this step is completed, one moves from the current case to another case according to the mechanism
described in points (1), (2), (3) and (4) in the proof. Then everything is repeated; this continues for a countable number of steps.

Now if we were to consider  the reverse of the permutation instead, then the description in the previous paragraph concerning  the contributions of the four cases  to the limiting object still  holds, but with 
the following transformation on the four cases: 
$$
\begin{aligned}
&\text{original}\to\text{reverse}:\\ 
& (+F)\to(-I);\   (-F)\to(+I);\  (+I)\to(-F);\  (-I)\to(+F).
\end{aligned}
$$
Similarly if we were to consider the complement of the permutation, then the description above would hold with the following transformation:
$$
\begin{aligned}
&\text{original}\to\text{complement}:\\
&(+F)\to(-F);\ (-F)\to(+F);\  (+I)\to(-I);\ (-I)\to(+I).
\end{aligned}
$$
And if we were to consider the reverse-complement of the  permutation, then  
the description above would hold with the following transformation:
$$
\begin{aligned}
&\text{original}\to\text{reverse-complement}:\\
&(+F)\to(+I); \ (-F)\to(-I);\  (+I)\to(+F); \ (-I)\to(-F).
\end{aligned}
$$

\begin{figure}\label{figure}
\begin{center}
\begin{tikzpicture}
\draw (0,0) to (0.2,0) to (0.2,0.2) to (1.2,0.2) to (1.2,1.2) to (0.2, 1.2) to (0.2,0.2) to (0,0.2) to  (0,0)
(2,1) to (3.2,1) to (3.2,0) to (2.2,0) to (2.2,1.2) to (2,1.2) to (2,1)
(4,0) to (5,0) to (5,1.2) to (5.2,1.2) to (5.2,1) to (4,1) to (4,0)
(6,.2) to (7.2, .2) to (7.2,0) to (7,0) to (7,1.2) to (6, 1.2) to (6,.2);
\end{tikzpicture}
\caption{$(+F)$, $(-F)$, $(+I)$, $(-I)$}
\end{center}
\end{figure}
Furthermore, we have the following connections between the limiting object obtained from  the original permutation
and the limiting objects obtained from the reverse, the complement and the reverse-complement of the permutation:

\noindent \it Reverse:\rm\
The   very block of numbers that
is contributed via  $(+F)$  to the limiting object 
for the original permutation is discarded via $(-I)$ and doesn't appear in the limiting object for the reverse of the permutation, and the very block of numbers  that is discarded via  $(-I)$ and doesn't appear in the limiting object for the original permutation is contributed via $(+F)$ to the limiting object for  the reverse
of the permutation. (This dictates  that the limiting random variable in Corollary \ref{cor} will be symmetric
with respect to $\frac12$.)

\noindent \it Complement:\rm\ The length of the block of numbers that is contributed via $(+F)$ to the limiting object 
for the original permutation is equal to the length of the block of infinities that is contributed via $(-F)$ to the limiting object  for the complement of the permutation, and the length of the block of infinities that is contributed
via $(-F)$ to the limiting object  for the original permutation is equal to the length of the block of numbers that is contributed via
$(+F)$ to the limiting object for the complement of the permutation.
 (This dictates that the limiting random variable in Corollary \ref{cortoprove} will be symmetric with respect to $\frac12$.)
 
\noindent \it Reverse-Complement:\rm\ The length of the block of infinities that is contributed via $(-F)$ to the  limiting object for the original permutation is equal to the length of the block of numbers that is discarded via $(-I)$ and won't appear in the limiting
object for the reverse-complement of the permutation, and the length of the block of numbers that is discarded via
$(-I)$ and won't appear in the limiting object for the original permutation is equal to the length of the block
of infinities that is contributed via $(-F)$ to the limiting object for the reverse-complement of the permutation.

\section{Proof of Proposition \ref{stable} and Corollary \ref{cortoprove} }\label{proppf}
\noindent \it Proof of Proposition \ref{stable}.
Part (i).\rm\ We need to show that
$E\exp(-it\frac{\sum_{k=1}^nL^{(k)}}{n^2})$ converges to the characteristic function appearing in part (i).
We consider $t$ to be fixed for the proof.
We have
\begin{equation}\label{charfuncsum}
E\exp(-it\frac{\sum_{k=1}^nL^{(k)}}{n^2})=\big(E\exp(-it\frac L{n^2})\big)^n.
\end{equation}
From \eqref{LR} and \eqref{genfunc} along with the fact that  \eqref{asymp} guarantees that \eqref{genfunc} holds with $x$ replaced by complex $z$ satisfying  $|z|=3-2\sqrt2$,
we have
\begin{equation}\label{charfunc}
\begin{aligned}
&E\exp(-\frac{i t}{n^2}L)=\frac{2(3-2\sqrt2)}{2-\sqrt2}\exp(-\frac{it}{n^2})+\sum_{j=2}^\infty
\frac{s_j(3-2\sqrt2)^j}{2-\sqrt2}\exp(-\frac{itj}{n^2})=\\
&\frac{(3-2\sqrt2)}{2-\sqrt2}\exp(-\frac{it}{n^2})+\frac1{2-\sqrt2}\thinspace s\big((3-2\sqrt2)e^{-\frac{it}{n^2}}\big)=
\frac{2-\sqrt2}2\exp(-\frac{it}{n^2})+\\
&\frac{2+\sqrt2}4\Big(1-(3-2\sqrt2)e^{-\frac{it}{n^2}}-
\sqrt{(3-2\sqrt2)^2e^{-\frac{2it}{n^2}}-6(3-2\sqrt2)e^{-\frac{it}{n^2}}+1}\ \Big).
\end{aligned}
\end{equation}
Noting that $3-2\sqrt2$ is a root of $z^2-6z+1$, we have after expanding in a power series and doing some arithmetic,
\begin{equation}\label{complexexpand}
(3-2\sqrt2)^2e^{-\frac{2it}{n^2}}-6(3-2\sqrt2)e^{-\frac{it}{n^2}}+1=(18\sqrt2-25)\frac{t^2}{n^4}+i(12\sqrt2-16)\frac t{n^2}+O(\frac1{n^6}).
\end{equation}
For   $n$ sufficiently large, this number lies in the right half plane. Thus we  interpret the square root
above as $\sqrt{z}=|z|^\frac12\exp(\frac12i \text{Arg}(z))$, with $\text{Arg}(z)\in(-\frac\pi2,\frac\pi2)$.
Thus, from \eqref{complexexpand} we have
$$
\sqrt{(3-2\sqrt2)^2e^{-\frac{2it}{n^2}}-6(3-2\sqrt2)e^{-\frac{it}{n^2}}+1}=
\frac{2(3\sqrt2-4)^\frac12|t|^\frac12}ne^{i\thinspace\text{sgn}(t)\frac\pi4}+O(\frac1{n^2}).
$$
Thus,
\begin{equation}
\begin{aligned}
&1-(3-2\sqrt2)e^{-\frac{it}{n^2}}-
\sqrt{(3-2\sqrt2)^2e^{-\frac{2it}{n^2}}-6(3-2\sqrt2)e^{-\frac{it}{n^2}}+1}=\\
&2\sqrt2-2-\frac{2(3\sqrt2-4)^\frac12|t|^\frac12}ne^{i\thinspace\text{sgn}(t)\frac\pi4}+O(\frac1{n^2}),
\end{aligned}
\end{equation}
and from \eqref{charfunc},
\begin{equation}\label{charfuncexpl}
E\exp(-\frac{i t}{n^2}L)=1-(\frac12)^\frac14\frac{|t|^\frac12}ne^{i\thinspace\text{sgn}(t)\frac\pi4}+O(\frac1{n^2}),
\end{equation}
(where we've used the fact that $\frac{2+\sqrt2}2(3\sqrt2-4)^\frac12=(\frac12)^\frac14$).
From \eqref{charfuncsum} and \eqref{charfuncexpl}, we obtain
\begin{equation}
\begin{aligned}
&\lim_{n\to\infty}E\exp(-it\frac{\sum_{k=1}^nL^{(k)}}{n^2})=\lim_{n\to\infty}\Big( 1-(\frac12)^\frac14\frac{|t|^\frac12}ne^{i\thinspace\text{sgn}(t)\frac\pi4}+O(\frac1{n^2})\Big)^{\frac1n}=\\
&\exp\big(-(\frac12)^\frac14|t|^{\frac12}e^{i\thinspace\text{sgn}(t)\frac\pi4}\big)=\exp\Big(-(\frac12)^\frac34|t|^\frac12\big(1+i\thinspace\text{sgn}(t)\big)\Big).
\end{aligned}
\end{equation}
\it \noindent Part (ii).\rm\
The proof is very similar to that of part (i),
so we leave it to the reader. However, for use in part (iii), we note that similar to \eqref{charfuncexpl}, we have
\begin{equation}\label{withR}
E\exp(-\frac{i t}{n^2}R)=1-2^\frac14\frac{|t|^\frac12}ne^{i\thinspace\text{sgn}(t)\frac\pi4}+O(\frac1{n^2}).
\end{equation}
\noindent \it Part (iii).\rm\ We will prove the result concerning
$\frac1{n^2}\sum_{k=1}^n\sum_{m=1}^{N^{(k)}}R_m^{(k)}$.
The proof for $\frac1{n^2}\sum_{k=1}^n\chi_{0,1}^{(k))}L^{(k)}$ is very similar.
What we will show is that the characteristic functions of the random variables
$\frac1{n^2}\sum_{k=1}^n\sum_{m=1}^{N^{(k)}}R_m^{(k)}$ converge to $\phi_{Z_R}(tEN^{(1)})$, where $\phi_{Z_R}$ is as in part (ii).
As noted after \eqref{N}, $EN^{(1)}=(\frac12)^\frac32$.
For the case of $\frac1{n^2}\sum_{k=1}^n\chi_{0,1}^{(k))}L^{(k)}$, one shows that the characteristic functions of these random variables converge
to $\phi_{Z_L}(tE\chi_{0,1}^{(1)})=\phi_{Z_L}(\frac12 t)$.
It turns out that both $\phi_{Z_R}((\frac12)^{\frac32}t)$ and $\phi_{Z_L}(\frac12t)$ are equal to the characteristic function in part (iii).

We consider $t$ fixed for the proof.
We have
\begin{equation}\label{beginformula}
E\exp\big(-\frac{it}{n^2}\sum_{k=1}^n\sum_{m=1}^{N^{(k)}}R_m^{(k)}\big)=\Big(E\exp\big(-\frac{it}{n^2}\sum_{m=1}^{N^{(1)}}R_m^{(1)}\big)\Big)^n,
\end{equation}
and from \eqref{N} and conditioning,
\begin{equation}\label{withcond}
E\exp\big(-\frac{it}{n^2}\sum_{m=1}^{N^{(1)}}R_m^{(1)}\big)=\frac{\sqrt2}2+\sum_{j=1}^\infty\sqrt2(3-2\sqrt2)^j\Big(E\exp(-\frac{it}{n^2}R)\Big)^j.
\end{equation}
Using \eqref{withR} and noting that $(1-z)^j=1-jz+R_2(z)$, with $|R_2(z)|\le \frac{j(j-1)}2|z|^2$, we have
\begin{equation}\label{longcalc}
\begin{aligned}
&\frac{\sqrt2}2+\sum_{j=1}^\infty\sqrt2(3-2\sqrt2)^j\Big(E\exp(-\frac{it}{n^2}R)\Big)^j=\\
&\frac{\sqrt2}2+\sum_{j=1}^\infty\sqrt2(3-2\sqrt2)^j\Big(1-2^\frac14\frac{|t|^\frac12}ne^{i\thinspace\text{sgn}(t)\frac\pi4}+O(\frac1{n^2})\Big)^j=\\
&\frac{\sqrt2}2+\sum_{j=1}^\infty\sqrt2(3-2\sqrt2)^j\big(1-2^\frac14\frac{|t|^\frac12}ne^{i\thinspace\text{sgn}(t)\frac\pi4}\thinspace j\big)+O(\frac1{n^2})=\\
&1-2^\frac14\frac{|t|^\frac12}ne^{i\thinspace\text{sgn}(t)\frac\pi4}EN+O(\frac1{n^2}).
\end{aligned}
\end{equation}
From \eqref{beginformula}. \eqref{withcond}, \eqref{longcalc} and the fact that $EN=(\frac12)^\frac32$, we obtain
$$
\lim_{n\to\infty}E\exp\big(-\frac{it}{n^2}\sum_{k=1}^n\sum_{m=1}^{N^{(k)}}R_m^{(k)}\big)=\exp\Big(-(\frac12)^\frac74|t|^\frac12\big(1+i\thinspace\text{sgn}(t)\big)\Big).
$$

\hfill$\square$

\medskip

\noindent \it Proof of Corollary \ref{cortoprove}.\rm\ By part (iii) of Proposition \ref{stable} along with the fact that
$\chi_{0,1}^{(k)}$ and $1-\chi_{0,1}^{(k)}$ have the same distribution, it follows that
both $\frac1{n^2}\sum_{k=1}^n\chi_{0,1}^{(k)}L^{(k)}$ and
$\frac1{n^2}\sum_{k=1}^n(1-\chi_{0,1}^{(k)})L^{(k)}$ converge in distribution as $n\to\infty$ to the distribution of $Z$
given in part (iii) of that proposition. Thus,
it remains to show  that these two sums, which for fixed  $n$ are not independent, are in fact asymptotically independent.
That is, we need to show that
\begin{equation}\label{toshow}
\begin{aligned}
&\lim_{n\to\infty}E\exp\big(-it\frac{\sum_{k=1}^n\chi_{0,1}^{(k)}L^{(k)}}{n^2}-
is\frac{\sum_{k=1}^n(1-\chi_{0,1}^{(k)})L^{(k)}}{n^2}\big)=\\
&\exp\Big(-(\frac12)^\frac74|t|^\frac12\big(1+i\thinspace\text{\rm sgn}(t)\big)\Big)
\exp\Big(-(\frac12)^\frac74|s|^\frac12\big(1+i\thinspace\text{\rm sgn}(s)\big)\Big).
\end{aligned}
\end{equation}
We have
\begin{equation}\label{1toshow}
\begin{aligned}
&E\exp\big(-it\frac{\sum_{k=1}^n\chi_{0,1}^{(k)}L^{(k)}}{n^2}-
is\frac{\sum_{k=1}^n(1-\chi_{0,1}^{(k)})L^{(k)}}{n^2}\big)=\\
&\Big(E\exp\big(-\frac{it}{n^2}\chi_{0,1}^{(1)}L-\frac{is}{n^2}(1-\chi_{0,1}^{(1)})L\big)\Big)^n.
\end{aligned}
\end{equation}
Also,
\begin{equation}\label{2toshow}
\begin{aligned}
&E\exp\big(-\frac{it}{n^2}\chi_{0,1}^{(1)}L-\frac{is}{n^2}(1-\chi_{0,1}^{(1)})L\big)=\\
&\frac12E\exp\big(-\frac{it}{n^2}L\big)+\frac12E\exp\big(-\frac{is}{n^2}L\big).
\end{aligned}
\end{equation}
Using \eqref{1toshow}, \eqref{2toshow} and
\eqref{charfuncexpl}, one easily obtains \eqref{toshow}.\hfill $\square$.

\medskip

\noindent \bf Acknowledgement.\rm\ The author thanks one of the referees for suggesting the inclusion of a discussion   on  the four symmetries that naturally occur in uniformly random separable permutations.

\end{document}